\documentclass[11pt,oneside]{article}
\usepackage{amssymb}
\usepackage{amsmath}
\usepackage{latexsym}
\usepackage{amsfonts}
\usepackage{amscd}
\usepackage{amsthm}
\usepackage{fancyhdr}

\usepackage[headings]{fullpage}

\newcommand{\gt}[1]{\mathfrak{#1}}

\newcommand{\RR}{{\mathbb R}}
\newcommand{\CC}{{\mathbb C}}
\newcommand{\ZZ}{{\mathbb Z}}
\newcommand{\ZZh}{\ZZ+\tfrac{1}{2}}
\newcommand{\hZZ}{\tfrac{1}{2}\ZZ}

\newcommand{\FF}{{\mathbb F}}
\newcommand{\ii}{{\bf i}}
\newcommand{\lab}{{\langle}}    
\newcommand{\rab}{{\rangle}}    

\newcommand{\LL}{{\Lambda}}

\newcommand{\Aut}{\operatorname{Aut}}
\newcommand{\Id}{\operatorname{Id}}

\newcommand{\PSL}{\operatorname{\textsl{PSL}}}    
\newcommand{\Sp}{\operatorname{\textrm{Spin}}}    

\newcommand{\MM}{\mathbb{M}}    
\newcommand{\Co}{\operatorname{\textsl{Co}}}    
\newcommand{\Suz}{\operatorname{\textsl{Suz}}}
\newcommand{\Ru}{\operatorname{\textsl{Ru}}}

\newcommand{\Cl}{\operatorname{Cliff}}

\newcommand{\Cm}{\operatorname{CM}}

\newcommand{\aco}{{A_{\Co}}}

\newcommand{\aru}{A_{{\Ru}}}  

\newcommand{\asuz}{A_{{\Suz}}}    

\newcommand{\cgset}{\Omega}   

\newcommand{\cej}{\jmath}

\newcommand{\ceJ}{J}

\newcommand{\vac}{{\bf 1}}      
\newcommand{\cas}{{\bf \omega}} 
\newcommand{\scas}{{\bf \tau}}  

\newtheorem{thm}{Theorem}

\newtheorem{prop}[thm]{Proposition}

\theoremstyle{definition}
\newtheorem*{defn}{Definition}

\theoremstyle{remark}

\theoremstyle{plain}

\theoremstyle{plain}

\begin{document}

\title{
    \textsc{Vertex operators and sporadic groups}
          }

\author{John F. Duncan
          \footnote{
          Harvard University,
          Department of Mathematics,
          One Oxford Street,
          Cambridge, MA 02138,
          U.S.A.}
          {}\footnote{
          Email: {\tt duncan@math.harvard.edu};\;
          homepage: {\tt
          http://math.harvard.edu/\~{}jfd/}
               }
               }

\date{October 21, 2007}



\maketitle

\begin{abstract}
In the 1980's, the work of Frenkel, Lepowsky and Meurman, along with
that of Borcherds, culminated in the notion of vertex operator
algebra, and an example whose full symmetry group is the largest
sporadic simple group: the Monster. Thus it was shown that the
vertex operators of mathematical physics play a role in finite group
theory. In this article we describe an extension of this phenomenon
by introducing the notion of enhanced vertex operator algebra, and
constructing examples that realize other sporadic simple groups,
including one that is not involved in the Monster.
\end{abstract}

\section{Motivation}\label{sec:motivation}

We begin not with the problem that motivates the article, but with
motivation for the tools that will furnish the solution to this
problem. The tools we have in mind are called {\em vertex operator
algebras (VOAs)}; here follows one way to motivate the notion.

In mathematics there are various kinds of finite dimensional
algebras that have proven to be significant or interesting in some
respect. For example,
\begin{enumerate}
\item     semisimple Lie algebras (with invariant bilinear form)
\item     simple Jordan algebras (of type A, B, or C)
\item     the Chevalley algebra (see \cite{CheSpinors})
\item     the Griess algebra (see \cite{GriFG})
\end{enumerate}
The items of this list are very different from each other in terms
of their properties and structure theory. Perhaps the only thing
they have in common (as algebras) is finite dimensionality.

Nonetheless, it turns out that there is such a process called {\em
affinization} which associates a certain infinite dimensional
algebra structure (let's say {\em affine algebra}), to each finite
dimensional example in this list. Not only this, but the
corresponding affine algebra has a distinguished representation
(infinite dimensional) for which the action of the affine algebra
extends in a natural way to a new kind of algebra structure called
vertex operator algebra structure\footnote{The VOA structure
corresponding to a semisimple Lie algebra with invariant bilinear form is
obtained in \cite{FreZhuVOAAffLieVir} (see also \cite{FLM}). VOAs
corresponding to the simple Jordan algebras of types A, B, and C are
given in \cite{LamJordanAlgVOAs}. An affinization of the Chevalley
algebra is constructed in \cite{FFR}, and the VOA corresponding to
the Griess algebra is constructed in \cite{FLM}.}.

Thus we obtain objects of a common category for each distinct
example here, and the notion of vertex operator algebra (VOA)
furnishes a framework within which these distinct finite dimensional
algebra structures may be unified.

\section{VOAs}

Let us now present a definition of the notion of VOA. For our
purposes it is more natural to consider the larger category of super
vertex operator algebras (SVOAs).

An SVOA is a quadruple $(U,Y,\vac,\cas)$ where
\begin{itemize}
\item     $U=U_{\bar{0}}\oplus U_{\bar{1}}$ is a super vector space
over a field $\FF$ say. (We will take $\FF$ to be $\RR$ or $\CC$.)
\item     $Y$ is a map $U\otimes U\to U((z))$, so that the image
of the vector $u\otimes v$ under $Y$ is a Laurent series with
coefficients in $U$. This series is denoted
$Y(u,z)v=\sum_{n}u_{(n)}vz^{-n-1}$, and the operator $Y(u,z)$ is
called the {\em vertex operator} associated to $u$.
\item     $\vac\in U_{\bar{0}}$ is called the {\em vacuum
vector}, and is a kind of identity for $U$ in the sense that we
should have $Y(\vac,z)u=u$ and $Y(u,z)\vac|_{z=0}u=u$ for all $u\in
U$.
\item     $\cas\in U_{\bar{0}}$ is called the {\em conformal
element}, and is such that for $Y(\cas,z)=\sum L(n)z^{-n-2}$ the
operators $L(n)$ should satisfy the relations
\begin{gather}\label{eqn:virrelns}
     [L(m),L(n)]=(m-n)L(m+n)+\frac{m^3-m}{12}\,c\,\delta_{m+n,0}\Id
\end{gather}
for some scalar value $c\in\FF$. In other words, the Fourier modes
of $Y(\cas,z)$ should generate a representation of the Virasoro
algebra, with some central charge $c$.
\end{itemize}
The main axiom for the vertex operators is the {\em Jacobi identity}, which states that
\begin{itemize}
\item     for $\ZZ/2$--homogeneous $u,v\in U$ and arbitrary $a\in U$ we should
have
\begin{gather}\label{eqn:cousin}
    \begin{split}
        &
        z_0^{-1}\delta\left(\frac{z_1-z_2}{z_0}\right)
        Y(u,z_1)Y(v,z_2)\\
    &-(-1)^{|u||v|}
        z_0^{-1}\delta\left(\frac{z_2-z_1}{-z_0}\right)
        Y(v,z_2)Y(u,z_1)
        \\
    &=
        z_2^{-1}\delta\left(\frac{z_1-z_0}{z_2}\right)
        Y(Y(u,z_0)v,z_2)
    \end{split}
\end{gather}
where $|u|$ is $1$ or $0$ as $u$ is odd or even, and similarly for
$|v|$, and the expression $z_0^{-1}\delta((z_1-z_2)/z_0)$, for example, denotes the formal power series
\begin{gather}\label{eqn:series}
     z_0^{-1}\delta\left(\frac{z_1-z_2}{z_0}\right)
     =\sum_{\substack{n,k\in\ZZ,\\ k\geq 0}}(-1)^k\binom{n}{k}z_0^{-n-1}z_1^{n-k}z_2^k.
\end{gather}

\end{itemize}
%
From (\ref{eqn:cousin}) and from the properties of $\vac$ we see the
extent to which the triple $(U,Y,\vac)$ behaves like a
(super)commutative associative (super)algebra with identity, since the formal series (\ref{eqn:series}) may be regarded as a ``delta function supported at $z_1-z_2=z_0$ (and expanded in $|z_1|>|z_2|$)''. The Jacobi identity (\ref{eqn:cousin}) thus encodes, among other things, some sense in which the compositions $Y(u,z_1)Y(v,z_2)$, $Y(v,z_2)Y(u,z_1)$, and $Y(Y(u,z_1-z_2)v,z_2)$ all coincide.


The structure in $(U,Y,\vac,\cas)$ which has no analogue in the
ordinary superalgebra case is that furnished by the conformal
element $\cas$. This structure furnished by $\cas$ (we will call it
{\em conformal structure}) manifests in two important axioms.
\begin{itemize}
\item     The action of $L(0)$ on $U$ should be diagonalizable with
eigenvalues in $\hZZ$ and bounded from below. We write
$U=\bigoplus_n U_n$ for the corresponding grading on $U$, and we
call $U_n$ the subspace of {\em degree} $n$.
\item     The operator $L(-1)$ should satisfy
$Y(L(-1)u,z)=D_zY(u,z)$ for all $u\in U$, where $D_z$ denotes
differentiation in $z$.
\end{itemize}
The conformal structure is essential for the construction of
characters associated to an SVOA. Zhu has shown \cite{ZhuPhd} that
these characters (under certain finiteness conditions on the SVOA)
span a representation of the Modular Group $\PSL_2(\ZZ)$.

We often write $U$ in place of $(U,Y,\vac,\cas)$. The scalar $c$ of
(\ref{eqn:virrelns}) is called the {\em rank} of $U$. Modules and
module morphisms can be defined in a natural way. We say that an
SVOA is {\em self-dual} in the case that it is irreducible as a
module over itself, and has no other inequivalent irreducible
modules.

\section{Sporadic groups} 

The Classification of Finite Simple Groups (see
\cite{GorLyoSolClassFSGs}, \cite{SolHstryClassFSG}) states that in
addition to the
\begin{itemize}
\item     cyclic groups of prime order
\item     alternating groups $A_n$ for $n\geq 5$
\item     finite groups of Lie type (like $\PSL_n(q)$, $PSU_n(q)$,
$G_2(q)$, \&c.)
\end{itemize}
there are exactly $26$ other groups that are finite and simple, and
can be included in none of the infinite families listed. These $26$
groups are called the {\em sporadic groups}.

The largest of the sporadic groups is called the {\em
Monster}\footnote{It is also called the {\em Friendly Giant}}, and
was first constructed by Robert L. Griess, Jr., who obtained this
result by explicitly constructing a certain commutative
non-associative algebra (with invariant bilinear form) of dimension
$196883$, with the Monster group $\MM$ as its full group of
automorphisms. This algebra is named the {\em Griess algebra}.

The Monster group furnishes a setting in which the majority (but not
all) of the other sporadic groups may be analyzed, since it involves
$19$ of the other sporadic simple groups. Here we say that a group
$G$ is {\em involved in the Monster} if $G$ is the homomorphic image
of some subgroup of the Monster; that is, if there is some $H$ in
$\MM$ with normal subgroup $N$ such that $H/N$ is isomorphic to $G$.
The sporadic groups that are involved in the Monster are called the
{\em Monstrous} sporadic groups.

The remaining $6$ sporadic groups not involved in the Monster are
called the {\em non-Monstrous} sporadic groups, or more colorfully,
the {\em Pariahs}.

\section{Vertex operators and the Monster} 

It is a surprising and fascinating fact that the Griess algebra can
appear in our list of finite dimensional algebras admitting
affinization (see \S\ref{sec:motivation}).

The fact that it does is due to the work of Frenkel, Lepowsky and
Meurmann, who extended the notion of affinization (known to exist
for certain Lie algebras) to the Griess algebra using vertex
operators, obtaining an infinite dimensional version: the affine
Griess algebra. They constructed a distinguished infinite
dimensional module over this affine Griess algebra called the
Moonshine Module. They built upon the work of Borcherds
\cite{BorPNAS} so as to arrive at the notion of vertex operator
algebra, and showed that the Moonshine Module admits such a
structure. Finally, they showed that the full automorphism group of
this structure is the Monster simple group.

In contrast to the situation with Lie algebras or Jordan algebras,
the Griess algebra is an object which is hard to axiomatize. It is
perhaps not clear that there is any reasonable category of algebras
(in the orthodox sense) which includes the Griess algebra as an
example (see \cite{ConCnstM}). The notion of VOA thus provides a
remedy to this situation: a setting within which the Griess algebra
can be axiomatized.

A related question is: ``How might the Monster group be
characterized?'' Having found that such an extraordinary group is
the symmetry group of some structure, we would like to be able
recognize this structure as distinguished in its own right, so that
our group might be defined to be just the group of automorphisms of
this distinguished object. In particular, our structure should
belong to some family of similar structures; a family equipped with
invariants, sufficiently rich that they can distinguish our
particularly interesting examples from all others, and sufficiently
simple that we can communicate them easily.

This question of characterization can also be addressed
(conjecturally, at least) within the theory of VOAs. Let us call to
attention three invariants for VOAs:
\begin{itemize}
\item     rank
\item     self-duality
\item     degree (vanishing conditions)
\end{itemize}
One may check that the Moonshine Module satisfies the following
three properties.
\begin{itemize}
\item     rank $24$
\item     self-dual
\item     degree $1$ subspace vanishes
\end{itemize}
It is a conjecture due to Frenkel, Lepowsky and Meurman that the
Moonshine Module is uniquely determined by these properties. Modulo
a proof of the conjecture, VOA theory thus provides a compelling
definition of the Monster group: the automorphism group of the
Moonshine Module, a beautifully characterized object in the category
of VOAs.

\section{Vertex operators and Monstrous groups}\label{sec:vopsmnstgps}

We have seen that vertex operators may be fruitfully applied to one
of the sporadic groups, and we may wonder if there is anything they
can do for the remaining $25$. Let us formulate The Problem:
\begin{itemize}
\item[$\ast$]     Given a sporadic group $G$, find a VOA whose
automorphism group is $G$, and characterize it.
\end{itemize}
The fact that $20$ of the sporadic groups are involved in the
Monster suggests that The Problem may have solutions, at least for
$G$ a Monstrous group. After all, if $G$ is such a group, then a
cover $\hat{G}=N.G$ say, of $G$, is a subgroup of the Monster, and
in particular, acts on the Moonshine Module. This action is probably
very reducible, but by choosing an appropriate irreducible
subalgebra for example, we may well obtain an object -- a VOA even
-- which serves as a reasonable analogue of the Moonshine Module for
our new group $\hat{G}$. We may even find that the normal subgroup
$N$ acts trivially, so that our analogue of the Moonshine Module
actually realizes $G$ itself, and not a cover of $G$. This outline
is extremely speculative, and certainly does not constitute an
acceptable solution to The Problem, but it does at least give us
somewhere to start, at least in the case that $G$ is a Monstrous
group.

It is important to mention that something like this has been carried
out rigourously and successfully in at least one case: that in which
$G$ is the Baby Monster $\textsl{BM}$, and the group $\hat{G}$ is a
double cover of $G$, the centralizer of a so-called $2A$ involution
in the Monster. The precise method is due to Gerald H\"ohn
\cite{HohnPhD}, and the result is a self-dual SVOA of rank
$23\tfrac{1}{2}$ whose full automorphism group is a direct product
$2\times \textsl{BM}$ of the Baby Monster with a group of order $2$.

\section{Vertex operators and the Conway group} 

At this point we would like to describe a solution to The Problem
for a specific Monstrous group $G$, which is nonetheless not along
the lines just described in \S\ref{sec:vopsmnstgps}. The group we
have in mind is the largest sporadic group of Conway, $\Co_1$. The
solution we have in mind is the object of the following Theorem.
\begin{thm}[\cite{DunVAco}]\label{thm:ACo_uniq}
Among nice rational $N=1$ SVOAs, there is a unique one satisfying
\begin{itemize}
\item     rank $12$
\item     self-dual
\item     degree $1/2$ subspace vanishes
\end{itemize}
\end{thm}
Let us name this structure $A_{\Co}$. We can see that it admits a
convenient characterization. That $A_{\Co}$ is a solution to one of
our problems is shown by the next Theorem.
\begin{thm}[\cite{DunVAco}]\label{thm:ACo_symm}
The full automorphism group of $A_{\Co}$ is the sporadic group
$\Co_1$.
\end{thm}
We should go no further before addressing the new terminology that
has arisen. The terms {\em nice} and {\em rational} refer to certain
technical conditions on SVOAs, and it will be convenient to put
aside their precise meaning, and refer the interested reader to the
article \cite{DunVAco}. (One would expect such technical conditions
to arise in a precise formulation of the uniqueness conjecture for
the Moonshine Module.) Of more importance for our present purpose is
the term {\em $N=1$ SVOA}.
\begin{defn}
An $N=1$ SVOA is a quadruple $(U,Y,\vac,\{\cas,\scas\})$, such that
$(U,Y,\vac,\cas)$ is an SVOA, and $\tau$ is a distinguished vector
of degree $3/2$ satisfying $\scas_{(0)}\scas=2\cas$, and such that
the Fourier coefficients of $Y(\tau,z)$ generate a representation of
the Neveu--Schwarz Lie superalgebra on $U$.\footnote{The definition
of $N=1$ SVOA here is almost identical to that of {\em $N=1$
Neveu--Schwarz vertex operator superalgebra without odd formal
variables ($N=1$ NS-VOSA)} which was introduced earlier by Barron
\cite{BarNSVOSA}.}
\end{defn}
The Neveu--Schwarz superalgebra is a natural super-analogue of the
Virasoro algebra. It is also known as the $N=1$ Virasoro
superalgebra. Thus an $N=1$ SVOA is just like an ordinary SVOA
except there is some extra structure: the role of the Virasoro
algebra is now played by the $N=1$ Virasoro superalgebra.

Let us consider for a moment longer, the difference between SVOA
structure and $N=1$ SVOA structure. Looking back at Theorem
\ref{thm:ACo_uniq}, experts will notice that the conclusion remains
true if we drop the the ``$N=1$'' from ``$N=1$ SVOA'' in the
hypothesis. That is to say, there is indeed a unique self-dual SVOA
with rank $12$ that has vanishing degree $1/2$ subspace. In fact, it
is a reasonably familiar object as SVOAs go: it is the lattice SVOA
associated to the integral lattice $D_{12}^+$ (the unique self-dual
integral lattice of rank $12$ with no vectors of unit norm).

One may be surprised to see a sporadic group, or even a finite group
here, since the SVOA underlying $A_{\Co}$ has infinite automorphism
group. In fact, there is an action by $\Sp_{24}$ (faithful up to
some subgroup of order $2$). The crux of the matter is that
\begin{enumerate}
\item     for a suitably chosen vector in this
$\Sp_{24}$-module $A_{\Co}$ the fixing group is $\Co_1$,
\item     the precise choice is made for us by the $N=1$ Virasoro superalgebra.
\end{enumerate}

Let us also emphasize that the uniqueness result for $A_{\Co}$
furnishes a compelling definition of the Conway group: as the full
automorphism group of $A_{\Co}$, a well characterized object in the
category of $N=1$ SVOAs.

\section{Enhanced SVOAs} 

With the example of $\aco$ in mind, and also with the suspicion that
it may be interesting to consider other extensions of the Virasoro
algebra, we formulate the notion of {\em enhanced SVOA}.

Roughly speaking, an enhanced SVOA is a quadruple
$(U,Y,\vac,\cgset)$ where $\cgset$ is a finite subset of $U$ (the
set of {\em conformal generators}) containing a vector $\cas$ for
which $(U,Y,\vac,\cas)$ is an SVOA. (We refer to
\cite[\S2]{DunVARuI} for the precise definition.) The subSVOA of $U$
generated by the elements of $\Omega$ is called the {\em conformal
subSVOA}.

The rank of an enhanced SVOA is just the rank of the underlying
SVOA. We say that an enhanced SVOA is self-dual just when it is
self-dual as an SVOA. The automorphism group of an enhanced SVOA is
the subgroup of the automorphism group of the underlying SVOA that
fixes every conformal generator.

We see then that an ordinary SVOA is an enhanced SVOA with
$\Omega=\{\cas\}$, and an $N=1$ SVOA is an enhanced SVOA with
$\Omega=\{\cas,\tau\}$, and conformal subSVOA a copy of the SVOA
associated to the vacuum representation of the $N=1$ Virasoro
superalgebra. In order to show that there are other interesting
examples of enhanced SVOA structure, we present the following
result.
\begin{thm}[\cite{DunVopsSpGps}]\label{thm:Asuz}
There exists a self-dual enhanced SVOA $\asuz$ of rank $12$
\begin{gather}
     \asuz=\left(\asuz,Y,\vac,\{\cas,\cej,\nu,\mu\}\right)
\end{gather}
with $\Aut(\asuz)\cong 3.\Suz$.
\end{thm}
It turns out that the conformal algebra in this example contains the
direct product of a pair of $N=1$ Virasoro superalgebras, at central
charges $11$ and $1$, respectively. The SVOA underlying $\asuz$
coincides with that underlying $\aco$, and taking the diagonal $N=1$
Virasoro superalgebra generated by $\tau=\nu+\mu$ (with central
charge $12$), we recover the enhanced SVOA structure with
automorphism group $\Co_1$.

\section{Beyond the Monster} 

We have seen already that vertex operators have a role to play in
the analysis of several Monstrous sporadic groups. In the cases of
Conway's group, Suzuki's group, the Monster, and the Baby Monster,
precise theorems have been formulated. A very significant question
is whether or not there is any application to sporadic groups beyond
the Monster; that is, to pariahs. We observed at the very beginning
that the notion of VOA can unify such disparate notions as
`semi-simple Lie algebra' and `commutative non-associative algebra'.

\medskip

The principal idea of this paper is that vertex operators do play a
role in the representation theory of non-Monstrous sporadic groups.

At this point let us introduce the sporadic group of Rudvalis,
$\Ru$. This sporadic simple group is not involved in the Monster; it
is one of the pariahs. It has order
\begin{gather}
     145926144000 = 2^{14}.3^3.5^3.7.13.29 \approx \frac{3}{2}\times
     10^{11}.
\end{gather}
The largest maximal subgroup is of the form $^2F_4(2)$ and has index
$4060$ in $\Ru$. (The Tits group has index two in this group.) The
next largest maximal subgroup is a non-split extension of the form
$2^6.G_2(2)$, and has index $188500$. The smallest non-trivial
irreducible representations of $\Ru$ have degree $378$.

\section{Vertex operators and Rudvalis's group} 

Consider The Problem for $G=\Ru$. The main theorem we wish to
present is the following.
\begin{thm}[\cite{DunVARuI},\cite{DunVopsSpGps}]\label{thm:Aru}
There exists a self-dual enhanced SVOA $\aru$ of rank $28$
\begin{gather}
     \aru=\left(\aru,Y,\vac,\{\cas,\cej,\nu,\varrho\}\right)
\end{gather}
with $\Aut(\aru)\cong 7\times \Ru$.
\end{thm}
(Compare this with the statement of Theorem \ref{thm:Asuz}.) We will
now provide a description of the enhanced SVOA $\aru$.

At the level of SVOAs, we have an isomorphism
\begin{gather}\label{eqn:arulattisom}
     \aru\cong V_{D_{28}^+}
\end{gather}
where $D_{28}^+$ denotes a self-dual integral lattice of rank $28$
with no vectors of unit length, and with $D_{28}$ as its even part.
We have observed already that there are analogous statements for the
SVOAs underlying the enhanced SVOAs $\aco$ and $\asuz$.
\begin{gather}
     \aco\cong\asuz\cong V_{D_{12}^+},
     \quad\text{ as SVOAs.}
\end{gather}
In the case of the enhanced SVOA $\asuz$, the conformal vectors
$\cej$ and $\cas$ have degree $1$ and $2$ respectively, and the two
conformal vectors beyond these; viz. $\nu$ and $\mu$, are both found
in the degree $3/2$ subspace. In the case of $\aru$, the degree
$3/2$ subspace is trivial. In fact the degree $1/2$, $3/2$, and
$5/2$ subspaces are all trivial for $\aru$. The extra conformal
vectors $\nu$ and $\varrho$ for $\aru$ are found in the degree $7/2$
subspace. This space is very large; the dimension is the number of
vectors of square-length $7$ in the lattice $D_{28}^+$.
\begin{gather}
     \dim (\aru)_{7/2}=2^{28}/2\approx 10^8
\end{gather}
The most effort in the construction of $\aru$ goes into determining
a precise description of the vectors $\nu$ and $\varrho$. It is a
remarkable fact that the finite group eventually obtained has almost
no other point-wise invariants\footnote{In fact, there is just one
other invariant in addition to $\nu$ and $\varrho$. It turns out to
be $\cej_{(0)}\nu$.} in its action on this $2^{27}$ dimensional
space.

\section{The $28$ dimensional representation} 

It is important for our construction of $A_{\Ru}$ that the Rudvalis
group admits a perfect double cover $2.\!\Ru$ which has irreducible
representations of degree $28$ (writable over $\ZZ[\ii]$). That the
group $2.\!\Ru$ preserves a lattice of rank $28$ over $\ZZ[\ii]$ was
observed independently by Meurman\footnote{private communication},
and by Conway and Wales \cite{ConWalRu} (see also \cite{ConRu} and
\cite{WilRu}), and this lattice is in fact self-dual when regarded
as a lattice (of rank $56$) over $\ZZ$. We choose to view this $28$
dimensional representation in terms of the maximal of $G_2$-type:
$2^6.G_2(2)$, which becomes $2^7.G_2(2)$ in the double cover
$2.\!\Ru$.

The action of the group $2^7.G_2(2)$ in the $28$ dimensional
representation can be understood in the following way in terms of
the $E_8$ lattice. Let $\LL$ denote a copy of the $E_8$ lattice, the
unique self-dual even lattice of rank $8$. (A {\em lattice} is a
free $\ZZ$-module equipped with a bilinear form, and an {\em even
lattice} is a lattice for which the square-norm of every vector is
an even integer.) We may take
\begin{gather}
     \LL=\left\{\sum_{i\in\Pi} n_ih_i\mid \sum_{i\in\Pi}
          n_i\in 2\ZZ;\;
          \text{all $n_i\in\ZZ$, or all $n_i\in\ZZh$}\right\}
\end{gather}
where the bilinear form is defined so that $\lab
h_i,h_j\rab=\delta_{ij}$. Then $\LL$ supports a structure of
non-associative algebra over $\ZZ$ which makes it a copy of the
integral Cayley algebra, or what is the same, a maximal integral
order in the Octonions. To see such a structure explicitly, we
assume that the index set $\Pi=\{\infty, 0,1,2,3,4,5,6\}$ is a copy
of the projective line over $\FF_7$, and then offer the following
defining relations (taken from \cite{ATLAS})
\begin{gather}
     {\bf 1}=\tfrac{1}{2}\sum_{i\in\Pi} h_i\\
     2h_{i}^2=h_i-{\bf 1}\\
     2h_{\infty}h_0={\bf 1}-h_3-h_5-h_6\\
     2h_0h_{\infty}={\bf 1}-h_2-h_1-h_4
\end{gather}
and also the images of these relations under the natural action of
$L_2(7)$ on the indices. The sublattice of doubles $2\LL$ is an
ideal in this algebra, and we may consider the quotient
$\bar{\LL}=\LL/2\LL$ which becomes a copy of the Cayley algebra over
the finite field $\FF_2$; what we call the {\em binary Cayley
algebra}. We should note that the automorphism group of this algebra
is the finite group $G_2(2)$ (which contains the simple group
$U_3(3)$ to index $2$).

There are just $2^8=256$ elements in the binary Cayley algebra. We
can count them.
\begin{gather}
\begin{tabular}{l|c}
  Type & Count \\\hline
  zeroes  & 1 \\
  identities & 1 \\
  involutions & 63 \\
  square roots of $0$ & 63 \\
  idempotents & 72 \\
  cube roots of $\bf{1}$ & 56
\end{tabular}
\end{gather}
There is a natural pairing on the elements of $\bar{\LL}$ obtained
by sending $x$ to the pair $\{x,{\bf 1}+x\}$. This association pairs
the zero with the identity, the involutions with the idempotents,
and partitions the idempotents and cube roots into $36$ and $28$
pairs, respectively. The group $G_2(2)$ acts transitively on each of
these different sets of pairs.

The $28$ cube root pairs are in a sense the basis upon which the
enhanced SVOA will be constructed. Let us denote them by $\Delta$. A
typical cube root of unity in $\bar{\LL}$ is $(h_i-h_j)$ for $i\neq
j\in\Pi$. A typical involution in $\bar{\LL}$ is given by
$(h_i+h_j)$ for $i\neq j\in\Pi$, and it is important that for any
given involution there are exactly $24$ cube roots ($12$ pairs of
cube roots) that are not orthogonal to the chosen involution. The
corresponding $12$-subsets of $\Delta$ are called {\em dozens}.

We now introduce a complex vector space $\gt{r}$ of dimension $28$,
with Hermitian form $(\cdot\,,\cdot)$ and an orthonormal basis
$\{a_i\}_{i\in\Delta}$, indexed by our cube root pairs. The
structures arising from the binary Cayley algebra described above
allow us to define an action by the group $2.2^6.G_2(2)$ on this
space. For example, the normal subgroup $2.2^6$ is generated by the
transformations that change sign on the coordinates of a given
dozen. With a somewhat finer analysis of the geometry of $\bar{\LL}$
we can define the action of the rest of the group (cf.
\cite{DunVARuI}); we call this group the {\em monomial group} and we
denote it $M$. (Warning: the action of $G_2(2)$ cannot be realized
as coordinate permutations. Instead we must write generators as
coordinate permutations followed by multiplications by $\pm 1$ or
$\pm \ii$ on particular coordinates. The group $M$ is a non-split
extension of $G_2(2)$.) The action of $M$ on $\gt{r}$ preserves the
Hermitian form.

\section{The conformal elements} 

Assume now that we have a Hermitian space $\gt{r}$ and a unitary
action on this space of the monomial group $M$ of the shape
$2.2^6.G_2(2)$. We assume further that $M$ is regarded as a matrix
group with respect to the basis $\{a_i\}$, and consists of monomial
matrices (having one non-zero entry in each row and column). We set
$\gt{u}=\gt{r}\oplus\gt{r}^*$, where $\gt{r}^*$ denotes the dual
space to $\gt{r}$, and is equipped with the induced Hermitian form.
The space $\gt{u}$ then comes equipped with a Hermitian form
(obtained by taking direct sum of those associated to the summands)
and also a bilinear form $\lab\cdot\,,\cdot\rab$, induced by the
canonical pairing $\gt{r}\times \gt{r}^*\to\CC$. We define
$\Cl(\gt{u})$ to be the Clifford algebra of $\gt{u}$ defined with
respect to this bilinear form.
\begin{gather}
     \Cl(\gt{u})=T(\gt{u})/\lab u\otimes u+\lab u,u\rab
          \mid u\in \gt{u}\rab
\end{gather}
We define $\Cm_X$ to be the module over $\Cl(\gt{u})$ spanned by a
vector $1_X$ satisfying $u1_X=0$ whenever $u\in\gt{r}^*$. We claim
that the isomorphism
\begin{gather}
     \Cm_X\cong \bigoplus \wedge^n(\gt{r})1_X
\end{gather}
holds when these spaces are viewed as modules over $\Cl(\gt{r})$
(the subalgebra of $\Cl(\gt{u})$ generated by
$\gt{r}\hookrightarrow\Cl(\gt{u})$). Next we claim that the degree
$7/2$ subspace of $\aru$ 
may be naturally identified with the even part of $\Cm_X$.
\begin{gather}
     (\aru)_{7/2}
          \longleftrightarrow
          \Cm_X^0\cong
          \bigoplus\wedge^{2n}(\gt{r})
\end{gather}
The Clifford algebra $\Cl(\gt{u})$ naturally contains a copy of the
group $\Sp(\gt{u})$, and the space $\Cm_X^0$ is an irreducible
module for this group $\Sp(\gt{u})$.

Recall that our goal is to define the elements $\nu$ and $\varrho$.
We now define $\nu$ by setting
\begin{gather}
     \nu=1_X+a_{\Delta}1_X
\end{gather}
where $a_{\Delta}=a_{\infty}a_1\cdots a_{27}\in\Cl(\gt{u})$ say. Let
us write $\overline{M}$ for the copy of $G_2(2)$ obtained by
replacing each non-zero entry with a $1$ in each matrix in $M$. The
vector $\varrho$ will be expressed in terms of orbits of
$\overline{M}$ on monomials $a_I1_X$ for $I\subset\Delta$ with
$|I|=14$.

It turns out that there are $80$ such orbits, but only $68$ give
rise to invariants for the monomial group $M$ in
$\wedge^{14}(\gt{r})1_X\subset\Cm_X^0$. The vector $\varrho$ is a
linear sum of these invariants, $t_i$ say, where the coefficients
$r_i$ can be taken to lie in $\ZZ[\ii]$.
\begin{gather}
     \varrho=\sum_{i=1}^{68} r_it_i
\end{gather}
The orbits are paired under complementation, and the coefficients of
invariants corresponding to complementary orbits are conjugate (up
to sign), so that ultimately, we require to specify $34$ values. We
refer to \cite[\S5.2]{DunVARuI} for the details.

Finally, we take $\aru=(\aru,Y,\vac,\{\cas,\cej,\nu,\varrho\})$ as
in the statement of Theorem \ref{thm:Aru}, and this completes the
construction.

It follows from the isomorphism (\ref{eqn:arulattisom}), with the
lattice SVOA for $D_{28}^+$, that $\aru$ is self-dual, and has rank
$28$. To prove that the automorphism group is of the stated form we
prove first that it is finite, by showing that it is a reductive
algebraic group with trivial Lie algebra (cf.
\cite[\S5.4]{DunVARuI}). It follows that it has dimension $0$, and
hence, is finite. We can explicitly construct generators for
$2.\!\Ru$ acting on $\aru$ and fixing the vectors
$\{\cas,\cej,\nu,\varrho\}$, and we may employ an argument from
\cite{NebRaiSloInvtsCliffGps}, to show that the only other
symmetries possible are scalar multiples of the identity. It is then
easy to check that such multiples must be $14^{{\rm th}}$-roots of
unity (since they must preserve $\varrho\in\wedge^{14}(\gt{r})1_X$).
Finally we obtain a central product $14\circ 2.\!\Ru$, but the
central $\ZZ/2$ here acts trivially on $\aru$, and so the full
automorphism group is just $7\times \Ru$.

\section{The character}

We conclude with consideration of the character of the enhanced SVOA
$\aru$.

The action of the Rudvalis group $\Ru$ on $\aru$ preserves a certain
vector $\cej$ of degree $1$. The residue of the corresponding vertex
operator $Y(\cej,z)$ is denoted $\ceJ(0)$, commutes with the
Virasoro operator $L(0)$, and has diagonalizable action on $\aru$,
thus giving rise to a grading by {\em charge}. It is natural then to
consider the two variable series
\begin{gather}
     {\sf tr}|_{\aru}p^{J(0)}q^{L(0)-c/24}
\end{gather}
which we call the {\em ($2$ variable) character} of $\aru$. Recall
the Jacboi theta function given by
\begin{gather}
     \vartheta_3(z|\tau)=\sum_{m\in\ZZ}e^{2\ii zm +\pi\ii\tau
     m^2}
\end{gather}
and also the Dedekind eta function
\begin{gather}
     \eta(\tau)=q^{1/24}\prod_{m\geq 1}(1-q^{m})
\end{gather}
written here according to the convention $q=e^{2\pi\ii\tau}$. Let us
also convene to write $p=e^{2\pi\ii z}$. Then we have
\begin{prop}[\cite{DunVopsSpGps}]
The character of $\aru$ is given by
\begin{gather}
     \begin{split}
     {\sf tr}|_{\aru}&p^{J(0)}q^{L(0)-c/24}=\frac{1}{2}\left(
          \frac{\vartheta_3(\pi z|\tau)^{28}}
          {\eta(\tau)^{28}}
          +\frac{\vartheta_3(\pi z+\pi/2|\tau)^{28}}
          {\eta(\tau)^{28}}\right)\\
          &+\frac{1}{2}p^{14}q^{7/2}\left(
          \frac{\vartheta_3(\pi z+\pi\tau/2|\tau)^{28}}
          {\eta(\tau)^{28}}
          +\frac{\vartheta_3(\pi z+\pi\tau/2+\pi/2|\tau)^{28}}
          {\eta(\tau)^{28}}\right)
     \end{split}
\end{gather}
\end{prop}
The terms of lowest charge and degree in the character of $\aru$ are
recorded in Table \ref{tab:RuMTseries}.
\begin{table}[t]
  \centering
  \caption{The Character of $\aru$}
  \label{tab:RuMTseries}
  \begin{small}
  \begin{tabular}{c|lllll}
     & 0 & 2 & 4 & 6 & 8\\\hline
    0  & 1 &   &   &  & \\
    1/2&   &   &   &   &\\
    1  & 784 & $378$ &  & & \\
    3/2 &  &  &  & & \\
    2 & 144452 & 92512 & 20475 & & \\
    5/2 &  &  &  & & \\
    3 & 11327232 & 8128792 & 2843568 & 376740 & \\
    7/2 & 40116600 & 30421755 & 13123110 & 3108105 & 376740 \\
    4 &  490068257 & 373673216 & 161446572 & 35904960 & 3108105\\
    9/2 & 2096760960 & 1649657520 & 794670240
     & 226546320 & 35904960\\
     5 & 13668945136 & 10818453324 & 5284484352 & 1513872360 & 226546320\\
     11/2 & 56547022140 & 45624923820 & 23757475560 & 7766243940 & 1513872360 \\
  \end{tabular}
  \end{small}
\end{table}
The column headed $m$ is the coefficient of $p^m$ (as a series in
$q$), and the row headed $n$ is the coefficient of $q^{n-c/24}$ (as
a series in $p$). The coefficients of $p^{-m}$ and $p^m$ coincide,
and all subspaces of odd charge vanish.

Many irreducible representations of $\Ru$ are visible in the entries
of Table \ref{tab:RuMTseries}. For example, we have the following
equalities, where the left hand sides are the dimensions of
homogeneous subspaces of $\aru$, and the right hand sides indicate
decompositions into irreducibles for the Rudvalis group.
\begin{gather}
     \begin{split}
     378&=378\\
     784&=1+783\\
     20475&=20475\\
     92512&=(2)378+406+91350\\
     144452&=(3)1+(3)783+65975+76125\\
     376740&=27405+65975+75400+102400
     \end{split}
\end{gather}

\newcommand{\etalchar}[1]{$^{#1}$}


\begin{thebibliography}{CCN{\etalchar{+}}85}

\bibitem[Bar00]{BarNSVOSA}
Katrina Barron.
\newblock {$N=1$} {N}eveu-{S}chwarz vertex operator superalgebras over
  {G}rassmann algebras and with odd formal variables.
\newblock In {\em Representations and quantizations (Shanghai, 1998)}, pages
  9--35. China High. Educ. Press, Beijing, 2000.

\bibitem[Bor86]{BorPNAS}
Richard Borcherds.
\newblock Vertex algebras, {Kac}-{Moody} algebras, and the {Monster}.
\newblock {\em Proceedings of the National Academy of Sciences, U.S.A.},
  83(10):3068--3071, 1986.

\bibitem[CCN{\etalchar{+}}85]{ATLAS}
J.~H. Conway, R.~T. Curtis, S.~P. Norton, R.~A. Parker, and R.~A. Wilson.
\newblock {\em Atlas of finite groups}.
\newblock Oxford University Press, Eynsham, 1985.
\newblock Maximal subgroups and ordinary characters for simple groups, With
  computational assistance from J. G. Thackray.

\bibitem[Che54]{CheSpinors}
Claude~C. Chevalley.
\newblock {\em The algebraic theory of spinors}.
\newblock Columbia University Press, New York, 1954.

\bibitem[Con77]{ConRu}
J.~H. Conway.
\newblock A quaternionic construction for the {R}udvalis group.
\newblock In {\em Topics in group theory and computation (Proc. Summer School,
  University Coll., Galway, 1973)}, pages 69--81. Academic Press, London, 1977.

\bibitem[Con85]{ConCnstM}
J.~H. Conway.
\newblock A simple construction for the {F}ischer-{G}riess monster group.
\newblock {\em Invent. Math.}, 79(3):513--540, 1985.

\bibitem[CW73]{ConWalRu}
J.~H. Conway and D.~B. Wales.
\newblock Construction of the {R}udvalis group of order
  {$145,\,926,\,144,\,000$}.
\newblock {\em J. Algebra}, 27:538--548, 1973.

\bibitem[Dun06a]{DunVARuI}
{John F.} Duncan.
\newblock Moonshine for {R}udvalis's sporadic group {I}.
\newblock arXiv:math.RT/0609449, September 2006.

\bibitem[Dun06b]{DunVopsSpGps}
John~F. Duncan.
\newblock {\em Vertex operators, and three sporadic groups}.
\newblock PhD thesis, Yale University, 2006.

\bibitem[Dun07]{DunVAco}
John~F. Duncan.
\newblock Super-moonshine for {C}onway's largest sporadic simple group.
\newblock {\em Duke Math. J.}, 139(2):255--315, 2007.

\bibitem[FFR91]{FFR}
Alex~J. Feingold, Igor~B. Frenkel, and John F.~X. Ries.
\newblock {\em Spinor construction of vertex operator algebras, triality, and
  {$E\sp {(1)}\sb 8$}}, volume 121 of {\em Contemporary Mathematics}.
\newblock American Mathematical Society, Providence, RI, 1991.

\bibitem[FLM88]{FLM}
Igor Frenkel, James Lepowsky, and Arne Meurman.
\newblock {\em Vertex operator algebras and the {M}onster}, volume 134 of {\em
  Pure and Applied Mathematics}.
\newblock Academic Press Inc., Boston, MA, 1988.

\bibitem[FZ92]{FreZhuVOAAffLieVir}
Igor~B. Frenkel and Yongchang Zhu.
\newblock Vertex operator algebras associated to representations of affine and
  {V}irasoro algebras.
\newblock {\em Duke Math. J.}, 66(1):123--168, 1992.

\bibitem[GLS94]{GorLyoSolClassFSGs}
Daniel Gorenstein, Richard Lyons, and Ronald Solomon.
\newblock {\em The classification of the finite simple groups}, volume~40 of
  {\em Mathematical Surveys and Monographs}.
\newblock American Mathematical Society, Providence, RI, 1994.

\bibitem[Gri82]{GriFG}
Robert~L. Griess, Jr.
\newblock The friendly giant.
\newblock {\em Invent. Math.}, 69(1):1--102, 1982.

\bibitem[H{\"o}h96]{HohnPhD}
Gerald H{\"o}hn.
\newblock {\em Selbstduale {V}ertexoperatorsuperalgebren und das
  {B}abymonster}, volume 286 of {\em Bonner Mathematische Schriften [Bonn
  Mathematical Publications]}.
\newblock Universit\"at Bonn Mathematisches Institut, Bonn, 1996.
\newblock Dissertation, Rheinische Friedrich-Wilhelms-Universit\"at Bonn, Bonn,
  1995.

\bibitem[Lam99]{LamJordanAlgVOAs}
Ching~Hung Lam.
\newblock On {VOA} associated with special {J}ordan algebras.
\newblock {\em Comm. Algebra}, 27(4):1665--1681, 1999.

\bibitem[NRS01]{NebRaiSloInvtsCliffGps}
Gabriele Nebe, E.~M. Rains, and N.~J.~A. Sloane.
\newblock The invariants of the {C}lifford groups.
\newblock {\em Designs, Codes and Cryptography}, 24:99--122, 2001.

\bibitem[Sol01]{SolHstryClassFSG}
Ronald Solomon.
\newblock A brief history of the classification of the finite simple groups.
\newblock {\em Bull. Amer. Math. Soc. (N.S.)}, 38(3):315--352 (electronic),
  2001.

\bibitem[Wil84]{WilRu}
Robert~A. Wilson.
\newblock The geometry and maximal subgroups of the simple groups of {A}.
  {R}udvalis and {J}. {T}its.
\newblock {\em Proc. London Math. Soc. (3)}, 48(3):533--563, 1984.

\bibitem[Zhu90]{ZhuPhd}
Yongchang Zhu.
\newblock {\em Vertex Operator Algebras, Elliptic Functions, and Modular
  Forms}.
\newblock PhD thesis, Yale University, 1990.

\end{thebibliography}
\end{document}